\documentclass[a4paper]{amsart}

\usepackage{amsfonts,amssymb,amsthm,amsmath}
\usepackage[totalwidth=17cm,totalheight=22cm]{geometry}
\usepackage[T1]{fontenc}

\newtheorem{prop}{Proposition}[section]
\newtheorem{defi}[prop]{Definition}
\newtheorem{lemma}[prop]{Lemma}
\newtheorem{sublemma}[prop]{Sub-Lemma}
\newtheorem{thm}[prop]{Theorem}
\newtheorem{cor}[prop]{Corollary}

\newtheorem{remark}[prop]{Remark}

\newtheorem{q}[prop]{Question}

\newcommand{\HSt}{\tilde{HS}}
\newcommand{\be}{\begin{eqnarray}}
\newcommand{\ee}{\end{eqnarray}}
\newcommand{\C}{{\mathbb C}}
\newcommand{\R}{{\mathbb R}}

\newcommand{\N}{{\mathbb N}}

\begin{document}

\title{Volume maximization and the extended hyperbolic space}
\date{v1 --- Aug 2009}
\author{Feng Luo
~and Jean-Marc Schlenker
}
\address{F. L.: Department of Mathematics, Rutgers University, Piscataway, NJ 08854, USA.}
\address{J.-M. S.:
Institut de Math{\'e}matiques de Toulouse, UMR CNRS 5219,
Universit{\'e} Toulouse III,
31062 Toulouse cedex 9, France.
}
\thanks{J.-M.S. was partially supported by the ANR program Repsurf : ANR-06-BLAN-0311.
 F.L. was partially supported by the NSF: NSF-DMS0604352.}

\begin{abstract}
We consider a volume maximization program to construct hyperbolic structures on
triangulated 3-manifolds, for which previous progress has lead to consider angle
assignments which do not correspond to a hyperbolic metric on each simplex.
We show that critical points of the generalized volume are associated to
geometric structures modeled on the extended hyperbolic space -- the
natural extension of hyperbolic space by the de Sitter space -- except for the
degenerate case where all simplices are Euclidean in a generalized sense.

Those extended hyperbolic structures can realize geometrically a decomposition of
the manifold as connected sum, along embedded spheres (or projective planes)
which are totally geodesic, space-like surfaces in the de Sitter part of the
extended hyperbolic structure.
\end{abstract}

\maketitle

\section{Introduction}

\subsection{The Casson-Rivin program and its extension}

\medskip

There are several ways to construct hyperbolic metrics on an
ideally triangulated 3-manifold with torus boundary.  The most
prominent ones are Thurston's algebraic gluing equations and
Casson-Rivin's angle structure.   In the angle structure approach,
one first introduces the notion of angles on the triangulation and
use the Lobachevsky function to define the volume of a tetrahedron
with angle assignments. In \cite{luo-ajm}, Casson-Rivin's angle
structure program is extended to closed triangulated 3-manifold.
The volume of angle structure is defined
and the critical points of the volume is investigated in \cite{luo-ajm}.

The goal of this paper
is to investigate the relationship between the critical points of
the volume and a natural extension of the hyperbolic space by the de Sitter space.
It turns out many of the critical
points have geometric meaning in terms of geometric structures based on
this extended hyperbolic space.

\subsection{Volume maximization on triangulated $3$-manifolds}

\begin{defi}
Let $M$ be a closed 3-manifold, and let $T$ be a triangulation of
$M$.  Recall that the triangulation $T$ is obtained as follows.
Take a finite set of disjoint Euclidean tetrahedra $s_1, ..., s_N$
and identify all their codimension-1 faces in pairs by affine
homeomorphisms. The quotient space $M$ inherits a natural
triangulation, denoted by $T$. A {\it wedge} of $T$ is a couple
$(e,s)$, where $s$ is a simplex of $T$ and $e$ is an edge of $s$
in the unidentified space $\cup_{i=1}^N s_i$. The set of wedges of
$T$ will be denoted by $W(T)$, while the set of vertices, edges,
2-faces and 3-simplices of $T$ will be denoted by $V(T), E(T),
F(T)$ and $S(T)$ respectively.
\end{defi}

The main ingredient in Casson-Rivin's angle structure is based on
the observation that the vertex link of an ideal hyperbolic
tetrahedron is a Euclidean triangle.  In our extension, we observe
that the vertex link of a compact tetrahedron in the hyperbolic,
spherical or Euclidean 3-space is a \it spherical \rm triangle.
This prompts us to propose the following definition.

\begin{defi}
An {\it angle structure} on $T$ is a function
$\theta:W(T)\rightarrow (0,\pi)$ such that:
\begin{itemize}
\item for each edge $e$ of $T$, the sum of the $\theta(w)$ over
the wedges of the form $w=(e,s), s\in S(T)$, is equal to $2\pi$,
\item for each $s\in S(T)$ and each vertex $v$ of $s$,
$\theta(e_1,s)+\theta(e_2,s)+\theta(e_3,s) >\pi$,
where the $e_i$ are the 3 edges of $s$ adjacent to $v$,
\item for each $s\in S(T)$ and each vertex $v$ of $s$,
$\theta(e_1,s)+\theta(e_2,s)< \theta(e_3,s)+\pi$, and similarly
for permutations of $1,2,3$.
\end{itemize}
We denote by $AS(T)$ the space of angle assignments on $T$.
\end{defi}

The geometric meaning of the first condition is that, if the
geometric structures on the simplices can be glued along the faces
(this is discussed below) then the angles add up at the edges so
that there is no singularity there. The second and the third
conditions mean that for each simplex the link of each vertex is a
spherical triangle.

It is quite possible that, in some cases, $AS(T)$ might be empty, however all manifolds
do admit a triangulation for which $AS(T)$ is non-empty \cite{kitaev-luo}.

There is a well-defined manner, explained in \cite{luo-ajm} and recalled in
sections 4-5, to associate to an angle assignment $\theta\in AS(T)$ a number which
is, in a precise way, a generalized volume. This defines a function
$V:AS(T)\rightarrow \R$. This ``volume'' is defined in terms
of a natural extension of the Schl\"afli formula, so that it automatically
verifies this identity.

If $\theta_0$ is a critical point of $V$ in $AS(T)$
such that all the angles assigned to all simplices of $T$
are the dihedral angles of a hyperbolic simplex, then, thanks to the Schl\"afli
formula, the lengths of an edge in the wedges containing it match, so that
the faces of the simplices can be glued
isometrically and $\theta_0$ defines a hyperbolic metric on $M$.
One of the main points of this paper is that an extension of this statement
holds when $\theta_0$ does not assign to all simplices of $T$ the dihedral
angles of a hyperbolic simplex.

\subsection{Extended hyperbolic structures}

There is a rather natural extension of the hyperbolic space by the de Sitter
space, used for instance in \cite{shu,cpt} in a polyhedral context somewhat
reminiscent of the arguments followed here. We call $HS^3$ this extended
hyperbolic space, so that $HS^3$ contains an open subset isometric to the
3-dimensional hyperbolic space and another open subset isometric to the quotient
of the de Sitter space by the antipodal map.

Given a 3-dimensional manifold $M$, an {\it HS-structure} on $M$
is a geometric structure locally modelled on $HS^3$. This is
explained in section 2.

\subsection{The main result}

The main result presented here is that most critical points of the volume
$V$ on the interior of $AS(T)$ have a natural interpretation in terms of
HS-structures.

\begin{thm} \label{tm:main}
Let $\theta\in AS(T)$ be a critical points of $V$. Then one of the
following applies.
\begin{enumerate}
\item $\theta$ corresponds to a spherical metric on $M$ (for each simplex
of $T$, the angles are the angles of a spherical simplex, and those simplex
can be glued isometrically along their faces, yielding a spherical metric),
\item $\theta$ corresponds to an
HS-structure on $M$. If this extended hyperbolic structure has a de
Sitter part, then it contains a totally geodesic surface
homeomorphic to $S^2$ or $\R P^2$, which is normal surface in $T$,
\item all simplices in $T$ are either Euclidean or flipped Euclidean (see
below). The volume $V(\theta)$ is then a non-negative integer multiple of $2\pi^2$.
If at least one simplex is flipped, then $V>0$ and $(M,T)$ contains a normal
surface homeomorphic to a sphere or a projective plane.
\end{enumerate}
\end{thm}

The first case should be clear, and can happen only if $M$ is diffeomorphic either
to $S^3$ or to $\R P^3$. In the second case, the totally geodesic space-like
surfaces in the de Sitter parts of the HS-structure realize
geometrically a decomposition of $M$ in irreductible parts, each of which
carries a hyperbolic metric. The third case is quite special, presumably
it can occur only in very specific cases, see below.

\subsection{Geometric interpretation}

The first idea in Theorem \ref{tm:main} is that considering HS-structures radically
simplifies the way in which one can find hyperbolic structures on 3-manifolds by
a critical point argument. Indeed, if $M$ is reducible, there might still be
a critical point of $V$ on $AS(T)$, corresponding not to a hyperbolic metric (which is
impossible if $M$ is reducible) but to a HS-structure, with a de Sitter part
corresponding to each incompressible sphere in $M$.

There is however a limit to this argument, as it stands. Given a HS-structure $h$ on
$M$, the set of its hyperbolic points is a domain $N\subset M$ on which the restriction
of $h$ defines a complete hyperbolic metric. The de Sitter parts of $h$ are topologically
either products of $S^2\times \R$ or products of a projective plane by an interval. This
means that each connected component of the boundary of $N$ has to be either a sphere
or a projective plane, which is a very restrictive condition.


\subsection{Further possible extensions}

The construction of HS-structures associated to critical points of $V$ on $AS(T)$
suggests that a further extension of the space of angle assignments $AS(T)$ should be
possible, allowing for instance for angle assignments such that the sum of angles at a vertex of
a simplex is equal to or less than $\pi$. Such angles assignments would corresponds geometrically (at
critical points of $V$) to triangulations with at least one vertex in the de Sitter part of
the HS-structure obtained. This line of investigation is not pursued here.

\section{Extended hyperbolic structures}

\subsection{The extended hyperbolic space}

One way to define this extension
is to consider a strictly convex quadric $Q$ in the projective 3-space $\R P^3$.
Given two distinct points $x,y\in \R P^3\setminus Q$, let $D$ be the projective
line containing $x$ and $y$. If $D$ intersects $Q$, let $a,b$ be the points
between $D$ and $Q$, chosen so that $a,x,y,b$ occur in this cyclic order on $D$.
Then define the Hilbert distance between $x$ and $y$ as
$$ d_H(x,y) = \frac{1}{2} \log [a,b;x,y]~, $$
where $[\cdot,\cdot;\cdot,\cdot]$ denotes the cross-ratio. If $D$ does not intersect $Q$, use the
same formula with $a,b$ replaced by the complex intersections of the line
containing $x,y$ with $Q$. This defines a ``distance'' in which the ball bounded
by $Q$ can be interpreted as a projective model of $H^3$, while the outside it
a projective model of the quotient of the de Sitter space $dS^3$ by the antipodal
map. In particular, $d_H(x,y)$ can be:
\begin{itemize}
\item real and negative, if $x,y$ are in the ball bounded by $Q$, and this defines inside $Q$
a projective model of the hyperbolic 3-dimensional space (known as the Klein model).
\item real and positive, if $x,y$ are outside the ball bounded by $Q$ and the line
joining them intersects $Q$ in two points. This line is then
time-like in $dS^3$.
\item in $i(0,\pi)$, if $x,y$ are outside $Q$ and the line containing them does
not intersect $Q$, this line is then space-like in $dS^3$.
\item in $i\pi/2+\R$, if $x$ is inside the ball bounded by $Q$ and $y$ is outside it.
\item $0$, if the line joining $x$ and $y$ is tangent to $Q$. This line is then
a light-like in $dS^3$.
\end{itemize}

The same construction also works in dimension 2, yielding the extended hyperbolic
plane $HS^2$.

\subsection{The double cover}

It is sometimes helpful to consider the double cover $\HSt^3$ of
$HS^3$. It is diffeomorphic to $S^3$, and has two hyperbolic
components each isometric to $H^3$, and one de Sitter component
isometric to the full de Sitter space $dS^3$. The same works in
any dimension. $\HSt^3$ is composed of two copies of the hyperbolic 3-space,
and of a one copy of the whole de Sitter space.

\subsection{Extended hyperbolic structures}

An {\it HS-structure} on a closed 3-dimensional manifold $M$ is a geometric structure
locally modeled on $HS^3$, with transformation group $PSL(2,\C)$.

If $h$ is an HS-structure on $M$, the set of points of $M$ where $h$ is locally
hyperbolic is an open subset of $M$, which we denote by $M_H$, and the restriction
of $h$ to $M_H$ is a complete hyperbolic metric. Similary, the set of points of $M$
where $h$ is locally modeled on the de Sitter space is an open subset of $M$, for
which we use the notation $M_{dS}$. Then $M_H\cup M_{dS}$ is dense is $M$, and its
complement is a surface.

\section{Triangles}

\subsection{The cosine formula}

\begin{defi}
We call $\cosh$ the restriction of the function $x\mapsto (e^x+e^{-x})/2$
to $\R_{<0}\cup [0,i\pi]\cup (i\pi+\R_{>0})$.
\end{defi}

With this definition, $\cosh$ is a bijection from its domain of definition to
$\R$.

\begin{lemma} \label{lm:cos}
Let $\alpha_1, \alpha_2, \alpha_3$ be the (interior) angles of a hyperbolic
triangle, and let $a_1,a_2,a_3$ be the length of the opposite edges. Then
\begin{equation} \label{eq:cosh}
\cosh(a_1)=\frac{\cos(\alpha_1)+\cos(\alpha_2)\cos(\alpha_3)}{\sin(\alpha_2)\sin(\alpha_3)}~.
\end{equation}
\end{lemma}

\subsection{The triangle inequality}

Consider now a triangle  in $S^2$, of angles $\alpha_1, \alpha_2, \alpha_3$.
It is a simple exercise to check that those angles satisfy the following
equation:
\begin{equation*}
\alpha_2+\alpha_3<\alpha_1+\pi \tag{$TI_1$}
\end{equation*}
Using the exterior angles rather than the interior angles, this equation can
be written as
$$ (\pi-\alpha_1) < (\pi-\alpha_2)+(\pi-\alpha_3)~, $$
which is the triangle inequality for the dual triangle in the sphere.

\subsection{A classification of M\"obius triangles}

Following \cite{luo-tams}, we consider here a generalization of
the notion of spherical (or hyperbolic) triangle.

\begin{defi}
A {\it M\"obius triangle} is a triple
$(\alpha_1,\alpha_2,\alpha_3)\in (0,\pi)^3$. Given a M\"obius triangle, its
{\bf edge lengths} are the complex numbers $(a_1,a_2,a_3)$ defined
by Equation (\ref{eq:cosh}), with the definition of $\cosh$ given
above.
\end{defi}

The rationale for the terminology used here is that, for any
triple $(\alpha_1, \alpha_2,\alpha_3)\in (0,\pi)^3$, there exists
a triangle in the complex plane bounded by three circles, unique
up to M\"obius transformation, so that its inner angles are the
$\alpha_i$'s. The constructions used below are however based
mostly on the extended hyperbolic plane and on real projective
geometry, rather than complex projective geometry. However
sticking to the terms ``M\"obius triangle'' should be helpful to
the reader insofar as it is closer to the previous works on the
subjects, e.g. \cite{luo-ajm}, \cite{luo-tams}.

\begin{lemma} \label{lm:triangles}
Let $T=(\alpha_1,\alpha_2,\alpha_3)$ be a M\"obius triangle, let $s=\alpha_1+\alpha_2+\alpha_3$,
and let $a_i$ be the edge lengths of $T$.
Exactly one of the following five cases applies.
\begin{enumerate}
\item $T$ is {\bf spherical}: $s>\pi$, and the triangle
inequalities $(TI_1),(TI_2),(TI_3)$ hold. Then $a_1,a_2,a_3\in
i(0,\pi)$.
\item $T$ is {\bf hyperbolic}: $s<\pi$. Then the triangle inequalities
$(TI_1),(TI_2),(TI_3)$ hold, and $a_1,a_2,a_3<0$.
\item $T$ is {\bf Euclidean}: $s=\pi$. Then the triangle
inequalities $(TI_1),(TI_2),(TI_3)$ hold, and $a_1=a_2=a_3=0$.
\item $T$ is {\bf flipped hyperbolic}: $\alpha_2+\alpha_3>\alpha_1+\pi$ (or
similarly after a permutation of $1,2,3$). Then $a_1<0$ while
$a_2,a_3\in i\pi + \R_{>0}$.
\item $T$ is {\bf flipped Euclidean}:
$\alpha_2+\alpha_3=\alpha_1+\pi$ (or similarly after a permutation
of $1,2,3$). Then $a_2=a_3=i\pi$ and $a_1=0$.
\end{enumerate}
\end{lemma}

The proof (which is elementary) is based on two preliminary
statements. Let $\alpha_1,\alpha_2,\alpha_3\in (0,\pi)$, and let
$s=\alpha_1+\alpha_2+\alpha_3$.

\begin{sublemma} \label{slm:1}
\begin{itemize}
\item If $\alpha_2+\alpha_3<\pi$, then $\cosh(a_1)<1$ (resp. $>1$) if and only if
$s>\pi$ (resp. $s<\pi$).
\item If $\alpha_2+\alpha_3>\pi$, then $\cosh(a_1)<1$ (resp. $>1$) if and only if
$(TI_1)$ holds (resp. $\alpha_2+\alpha_3>\pi+\alpha_1$).
\end{itemize}
\end{sublemma}

\begin{sublemma} \label{slm:2}
$\cosh(a_1)>-1$ if and only if either $\alpha_2>\alpha_3$ and $(TI_3)$ holds,
or $\alpha_3>\alpha_2$ and $(TI_2)$ holds.
\end{sublemma}

\begin{proof}[Proof of Lemma \ref{lm:triangles}]
Suppose first that the three triangle inequalities hold. Then $T$ is in one of the
cases (1), (2) or (4) depending on $s$. Note also that if $(TI_1)$ does not hold,
then clearly both $(TI_2)$ and $(TI_3)$ are satisfied. So $T$ is in case
(3) if there is equality in inequality $(TI_1)$, or (5) otherwise.

In case (1), $\cosh(a_1)<1$ by Sub-Lemma \ref{slm:1}, and $\cosh(a_1)>-1$ by
Sub-Lemma \ref{slm:2}, so $a_1\in i(0,\pi)$ by the definition of $\cosh$ used
here. The same holds of course for $a_2$ and $a_3$. In case (2), Sub-Lemma
\ref{slm:1} shows that $\cosh(a_i)=1$ for $i=1,2,3$ so that all the $a_i$ are zero.
In case (4), the first case of Sub-Lemma \ref{slm:1}, $\cosh(a_1)>1$, so that
$a_1<0$ by our definition of $\cosh$, and the same applies to $a_2$ and $a_3$.

In case (3), the second point in Sub-Lemma \ref{slm:1} shows that $\cosh(a_1)=1$,
while Sub-Lemma \ref{slm:2} shows that $\cosh(a_2)=\cosh(a_3)=-1$, so that
$a_2=a_2=i\pi$. The same argument applies to case (5), then $\cosh(a_1)>1$
while $\cosh(a_2),\cosh(a_3)<-1$, so that $a_1<0$ while $a_2,a_3\in i\pi + \R_{>}$.
\end{proof}

The following lemma shows that the edge lengths determine the shape of a M\"obius
triangle.

\begin{lemma} \label{lm:lengths-unique}
Let $A=(\alpha_1,\alpha_2,\alpha_3)$ and $B=(\beta_1, \beta_2,
\beta_3)$ be two M\"obius triangles which are not Euclidean or
flipped Euclidean. If the corresponding edge lengths of $A$ and
$B$ are the same, then $\alpha_i=\beta_i$ for all i.
\end{lemma}

The proof follows from the cosine law that
\begin{equation} \label{eq:cos}
\cos(\alpha_1)=\frac{-\cosh(a_1)+\cosh(a_2)\cosh(a_3)}{\sinh(a_2)\sinh(a_3)}.
\end{equation}
This laws shows that lengths determine the angles.

\subsection{Geometric realization of triangles}

The classification of M\"obius triangles in Lemma \ref{lm:triangles} has a natural
interpretation in terms of triangles in the extended hyperbolic plane $\HSt^2$.
There is no interpretation needed for spherical, Euclidean or hyperbolic triangle,
but flipped hyperbolic triangle correspond to triangles in $\HSt^2$ with one vertex
and two in the other.

Suppose $t$ is such a flipped hyperideal triangle, with vertices $v_1,v_2,v_3$,
with $v_1$ in one copy of $H^2$ and $v_2,v_3$ in the other. Let $\alpha_i$ be
the angle of $t$ at $v_i$, $1\leq i\leq 3$. Those angles can be understood by
``flipping'' $t$, that is, considering the triangle $t'$ with vertices
$v'_1, v_2,v_3$, where $v'_1$ is the antipode of $v_1$ in $\HSt^2$. $t'$ is
a ``usual'' hyperbolic triangle, and its angles are $\beta_1=\alpha_1$,
$\beta_2=\pi-\alpha_2, \beta_3=\pi-\alpha_3$. Since $t'$ is a hyperbolic
triangle, its angles satisfy
$$ \beta_1+\beta_2+\beta_3<\pi~, $$
which translates as
$$ \alpha_1+\pi<\alpha_2+\alpha_3~, $$
the inverse triangle inequality for $t$. Similarly, $t$ satisfies the triangle
inequality,
$$ \pi+\beta_2 > \beta_1+\beta_3~~ \mbox{and}~~ \pi+\beta_3>\beta_1+\beta_2~, $$
which translates as
$$ \alpha_1+\alpha_2 < \pi+\alpha_3~~\mbox{and}~~\alpha_1+\alpha_3<\pi+\alpha_2~. $$
This shows that $(\alpha_1,\alpha_2,\alpha_3)$ is a flipped hyperbolic triangle,
according to Lemma \ref{lm:triangles}. The same argument can be used backwards,
to show that any flipped hyperbolic triangle in the sense of  Lemma \ref{lm:triangles}
is the triple of angles of a flipped hyperbolic triangle in $\HSt^2$.

Flipped Euclidean triangles
can be understood in the same way but by taking a limit. The usual Euclidean
triangles can be considered as limits of hyperbolic triangles of diameter
going to zero -- actually, a blow-up procedure is necessary, since what we
really want to consider are sequences of degenerating hyperbolic triangles
for which the angles, and therefore the ratio of the edge lengths, has
a limit. This can also be done for flipped hyperbolic triangles, with
vertices converging either to a point in the one copy of the hyperbolic
plane of $\HSt^2$ or to its antipode in the other copy.

Putting this together, we obtain the following statement.

\begin{lemma} \label{lm:mobius}
Let $T=(\alpha_1, \alpha_2, \alpha_3)$ be a M\"obius triangle.
\begin{enumerate}
\item If $T$ is spherical, there is a triangle $t\subset S^2$, unique up
to the action of $O(3)$, with angles $\alpha_1,\alpha_2$ and $\alpha_3$.
\item If $T$ is hyperbolic, there is a triangle $t\in H^2$, unique
up to the hyperbolic isometries, with angles $\alpha_1,\alpha_2, \alpha_3$.
\item If $T$ is Euclidean, there is a triangle $t\subset \R^2$, unique
up to isometries and homotheties, with angles $\alpha_1,\alpha_2, \alpha_3$.
In other terms, there is a sequence of hyperbolic triangles $(t_n)_{n\in \N}$
in $H^2$, with angles $\alpha_{1,n},\alpha_{2,n},\alpha_{3,n}$
converging to $\alpha_1,\alpha_2,\alpha_3$, respectively.
\item If $T$ is flipped hyperbolic, there is a $t\in \HSt^2$, with one vertex
in one copy of $H^2$ and two in the other, with angles
$(\alpha_1, \alpha_2, \alpha_3)$. It is unique up to the action of
$O(2,1)$ on $\HSt^2$.
\item If $T$ is flipped Euclidean, there is a sequence of flipped hyperbolic
triangles $(t_n)_{n\in \N}$ in $\HSt^2$, with angles
$\alpha_{1,n},\alpha_{2,n},\alpha_{3,n}$
converging to $\alpha_1,\alpha_2,\alpha_3$, respectively.
\end{enumerate}
\end{lemma}

\section{Three-dimensional simplices}

\subsection{Angle System }

We now consider 3-dimensional simplices, and assign to each wedge an angle,
as follows.

\begin{defi}
Suppose $s$ is a tetrahedron.  Then an \it angle system \rm on $s$
is a function $\theta: \{ (e, s) | e  \text{ is an edge of  } s\}
\to (0, \pi)$ so that for three edges $e_i, e_j, e_k$ ending at a
vertex $v$ of $s$, the 3 angles $\theta(e_i, s), \theta(e_j, s),
\theta(e_k, s)$ are the inner angles of a spherical triangle.
Let $AS(s)$ be the space of all angle systems on $s$.  An {\it angled
3-simplex} is a 3-simplex together with an angle system.

If $T$ is a triangulation of a closed 3-manifold $M$, an angle
system on $T$ is a function $\theta$ defined on the set of all
wedges $\{ (e, s) |  e$ is an edge of a tetrahedron $s \}$  so
that
\begin{itemize}
\item for each edge $e$ of $T$, the sum of the values of $\theta$
on the wedges having $e$ as their edge is equal to $2\pi$, \item
for each 3-simplex $s$ in $T$, the restriction of $\theta$ to all
wedges in $s$ forms an angle system.
\end{itemize}
\end{defi}

In the paper \cite{luo-ajm}, the geometric prototype of an angled
3-simplex is the M\"obius tetrahedron. Namely a  topological
tetrahedron in $\R^3$ bounded by four 2-spheres of inner angles
less than $\pi$.  However, there are angled 3-simplex which cannot
be realized as a M\"obius 3-simplex. Our main observation is that,
in terms of HS-geometry, these angled 3-simplices all have a geometric
meaning. Furthermore, the edge lengths, volume and Schl\"afli
formula can be generalized to the HS-geometry. These
generalizations are exactly the underlying geometric meaning of
the corresponding notions defined in \cite{luo-ajm}.

\subsection{Face angles}

\begin{defi}
Let $\alpha=(\alpha_{12},\cdots,\alpha_{34})\in AS(s)$. The face
angles of $\alpha$ are the numbers $\beta_{ijk}\in (0,\pi)$
defined, for $\{ i,j,k,l\}=\{ 1,2,3,4\}$, by the formula
$$ \cos(\beta_{jk}^i)=\frac{\cos(\alpha_{il})+\cos(\alpha_{ij})\cos(\alpha_{ik})}
{\sin(\alpha_{ij})\sin(\alpha_{ik})}~. $$
\end{defi}

The geometric meaning of the face angle is as follows. According
to the definition, at the $i$th vertex $v_i$,   the angles
$\alpha_{ij}, \alpha_{ik}, \alpha_{il}$ are the inner angles of a
spherical triangle $A_{jkl}$, which can be considered as the link
$v_i$ in the tetrahedron.  Then the face angle $\beta_{jk}^i$ is
the $jk$-th edge length of $A_{jkl}$. By definition, face angles
are then in $(0,\pi)$.

\subsection{Edge lengths}

Using the faces angles, we make each codimension-$1$ face of an
angled tetrahedron $s$ a M\"obius triangle.   Thus, by lemma 3.2,
we can define for each edge in each face, an edge length. The
following is proved in \cite{luo-ajm}.

\begin{lemma} \label{lm:lengths-same}
If $L$ is an edge of a tetrahedron $s$ with an angle system and
$D_1, D_2$ are two faces of $s$ having $L$ as an edge, then the
length of $L$ in $D_1$ is the same as that in $D_2$.
\end{lemma}

Thus the following is well defined.

\begin{defi}
Let $\alpha=(\alpha_{12},\cdots,\alpha_{34})\in AS(s)$. The edge
lengths of $\alpha$ are the numbers $(l_{ij})_{i\neq j}$ defined
as follows: $l_{ij}$ is the length of the edge $ij$ in the two
faces of $T$ adjacent to the vertices $i$ and $j$.
\end{defi}

\subsection{A classification of simplices}

It is now possible to describe a classification of 3-dimensional
angled simplices. It is slightly more elaborate than the
corresponding classification for M\"obius triangles, because
simplices can be ``flipped'' in two ways, depending on whether one
or two vertices are in one of the copies of $H^3$ in $\HSt^3$.
Here is the definition of the flip at the i-th vertex of $\alpha
\in AS(s)$. See \cite{luo-ajm} for more details.  The flipped
simplex $\alpha' = (\alpha'_{12}, ..., \alpha'_{34})$ has angles
$\alpha'_{ij} = \alpha_{ij}$ for $j \neq i$ and
$\alpha_{jk}'=\pi-\alpha_{jk}$ for $j, k \neq i$. Geometrically,
if $v_1, v_2, v_3, v_4$ are the vertices of a spherical simplex,
then the flipped simplex (about the first vertex $v_1$) is the
spherical simplex with vertices $-v_1, v_2, v_3, v_4$ where $-v_1$
is the antipodal point of $v_1$.

\begin{lemma}\label{lm:simplex}
Let $\alpha\in AS(s)$. After a permutation of $\{ 1,2,3,4\}$,
$\alpha$ is of exactly one of the following types:
\begin{enumerate}
\item {\bf spherical}: all faces of $T$ are spherical triangles, and all
edge lengths are in $i(0,\pi)$.
\item {\bf hyperbolic}: all faces of $T$ are hyperbolic triangles, and all
edge lengths are in $\R_{<0}$.
\item {\bf flipped hyperbolic}: the face $(234)$ is a hyperbolic triangle,
while the faces adjacent to the vertex $1$ are flipped hyperbolic triangles.
The lengths of the edges $(12),(13),(14)$ are in $i\pi+\R_{>0}$, while the
length of the other edges are in $\R_{<0}$.
\item {\bf doubly flipped hyperbolic}: all faces of $T$ are flipped hyperbolic
triangles, the lengths of the edges $(12)$ and $(34)$ are negative numbers
while the length of the other edges are in $i\pi+\R_{>0}$.
\item {\bf Euclidean}: all faces of $T$ are Euclidean triangles, all edge lengths
are zero.
\item {\bf flipped Euclidean}: the length of the edges $(12),(13),14)$ are equal
to $i\pi$, while the lengths of $(14),(24),(34)$ are zero.
\item {\bf doubly flipped Euclidean}: the lengths of $(12),(34)$ are zero while
all other edges have length $i\pi$.
\end{enumerate}
\end{lemma}

The terminology is based, as for triangles, on the idea of ``flipping'' a
hyperbolic simplex: this means replacing one of its vertices, say $v_1$,
by its antipode in $\HSt^3$. The dihedral angles at all three edges not
adjacent to $v_1$ are then replaced by their complement to $\pi$, and
the effect on the edge lengths is as described in Lemma \ref{lm:simplex}.
Doubly flipping a hyperideal simplex means replacing two vertices by
their antipodes in $\HSt^3$.

\begin{proof}[Proof of Lemma \ref{lm:simplex}]
Let $\alpha\in AS(s)$, let $v_1, v_2, v_3, v_4$ be the vertices of $s$.
We consider different cases, depending on the lengths of the edges of $s$
and in particular of its face $(v_1,v_2,v_3)$.
\begin{enumerate}
\item $(v_1,v_2,v_3)$ is a spherical triangle, i.e., its edge lengths
are in $i(0,\pi)$. Lemma \ref{lm:triangles}, applied to the three
other faces of $s$, shows that the lengths of the three other edges
of $s$ are also in $i(0,\pi)$, it follows that $s$ is spherical.
\item $(v_1,v_2,v_3)$ is hyperbolic. Then its edge lengths are negative,
and considering Lemma \ref{lm:triangles} shows that there are only
two possibilities.
\begin{enumerate}
\item $(v_2,v_3,v_4)$ is hyperbolic, that is, its edge lengths are
negative. Then, again by Lemma \ref{lm:triangles}, the length of the
edge $(v_1,v_4)$ is also negative, so that $s$ is hyperbolic.
\item $(v_2,v_3,v_4)$ is flipped hyperbolic, that is, the lengths of
the edges $(v_2,v_4)$ and $(v_3,v_4)$ are in $i\pi+\R_{>0}$. Then the length
of $(v_1,v_4)$ is also in $i\pi+\R_{>0}$, so that $s$ is flipped hyperbolic.
\end{enumerate}
\item $(v_1,v_2,v_3)$ is flipped hyperbolic, we can suppose without loss of
generality that the length of $(v_1,v_2)$ is in $\R_{>0}$ and the lengths
of the two other edges are in $i\pi+\R_{>0}$. Two cases are possible.
\begin{enumerate}
\item $(v_1,v_2,v_4)$ is hyperbolic, it then follows from Lemma \ref{lm:triangles}
that $s$ is flipped hyperbolic.
\item $(v_1,v_2,v_4)$ is flipped hyperbolic, it then follows from Lemma \ref{lm:triangles}
that $s$ is doubly flipped hyperbolic.
\end{enumerate}
\item  $(v_1,v_2,v_3)$ is Euclidean, so that all its edges have zero length.
Lemma \ref{lm:triangles} then shows that there are two possible cases.
\begin{enumerate}
\item $(v_1,v_2,v_4)$ is Euclidean. Then all edges of $s$ have zero length, and
$s$ is Euclidean.
\item $(v_1,v_2,v_4)$ is flipped Euclidean, so that $(v_1,v_4)$ and $(v_2,v_4)$
have length $i\pi$. In this case $s$ if flipped Euclidean.
\end{enumerate}
\item $(v_1,v_2,v_3)$ is flipped Euclidean, we can suppose without loss of
generality that $(v_1,v_2)$ has zero length zero while the lengths of the
other two edges are equal to $i\pi$. There are again two cases to consider.
\begin{enumerate}
\item $(v_1,v_2,v_4)$ is Euclidean, so that all its edge lengths are zero, and
it easily follows that $s$ is flipped Euclidean.
\item $(v_1,v_2,v_4)$ is flipped Euclidean, so that the edges $(v_1,v_4)$ and
$(v_2,v_4)$ have length $i\pi$. Then $s$ is doubly flipped Euclidean.
\end{enumerate}
\end{enumerate}
\end{proof}

According to the lemma, there are three types of angled
simplices: Euclidean, hyperbolic and spherical. A angled simplex
is of Euclidean (or hyperbolic) type if it can be flipped to a
Euclidean (or hyperbolic) simplex. A spherical type simplex is the
same as a spherical simplex.  The type of a simplex can be
determined by the length of any of its edges.

\begin{cor}
Suppose $e$ is an edge of an angled simplex of length $l(e)$. Then
$s$ is of
\begin{enumerate}
\item{\bf Euclidean type}   if and only if $l(e) \in \{ 0, i \pi
\}$,
\item{\bf hyperbolic type}  if and only if $l(e) \in \R_{<0}
\cup \{ i \pi + r | r \in \R_{>0}\}$,
\item{\bf spherical type} if and only if $l(e) \in  i(0, \pi)$.
\end{enumerate}
In particular, if $e$ and $e'$ are two edges of an angled simplex
$s$, then their lengths $l(e)$ and $l(e')$ are in the same subset
listed above.
\end{cor}

\subsection{Combinatorics of the space of simplices}

The classification given in \S4.4 can be interpreted in terms of
the HS-geometry as follows, as for M\"obius triangles in Lemma \ref{lm:mobius}.
Let $s$ be a simplex, and let $\alpha\in AS(s)$.
\begin{enumerate}
\item If $\alpha$ is spherical, the $\alpha_{ij}$ are the dihedral
angles of a unique simplex in $S^3$.
\item If $\alpha$ is hyperbolic,  the $\alpha_{ij}$ are the dihedral
angles of a unique simplex in $H^3$.
\item If $\alpha$ is hyperbolic, the $\alpha_{ij}$ are the dihedral
angles of a unique simplex in $\HSt^3$, with three vertices in one of
the copies of $H^3$ and one in the other.
\item  If $\alpha$ is hyperbolic, the $\alpha_{ij}$ are the dihedral
angles of a unique simplex in $\HSt^3$, with two vertices in one of
the copies of $H^3$ and two in the other.
\item If $\alpha$ is Euclidean, the $\alpha_{ij}$ are the dihedral
angles of an Euclidean simplex, unique up to homothety. They are
also limits of sequences of angles of hyperbolic simplices.
\item If $\alpha$ is flipped Euclidean, it is the limit of a
sequence of angles of flipped hyperbolic polyhedra.
\item Similarly, if $\alpha$ is doubly flipped Euclidean, it is the limit of a
sequence of angles of doubly flipped hyperbolic polyhedra.
\end{enumerate}

Consider now $AS(s)$ as the space of 6-tuples of angles in $(0,\pi)$
satisfying some linear inequalities. It contains some subdomains corresponding
to the different types of simplices. It is interesting to consider the
combinatorics of this decomposition of $AS(s)$. The definitions show clearly
that any continuous path going from a simplex of spherical type to a simplex
of hyperbolic type has to go through a simplex of Euclidean type. Moreover,
the only way to go from a hyperbolic simplex to a doubly hyperbolic simplex
is through spherical simplices, and similarly for doubly hyperbolic simplices.

\section{The generalized volume}

\subsection{The Schl\"afli formula}

The last part of the picture considered here is the generalized volume, which
is defined for the simplices in the extended hyperbolic space. There are
severaly ways to define it, we use here the Schl\"afli formula, which we
first recall for ``usual'' (spherical or hyperbolic) simplices. We refer
to \cite{milnor-schlafli,geo2} for a proof.

\begin{lemma} \label{lm:schlafli}
For any one-parameter family of spherical (resp. hyperbolic)
simplices, its volume $V$ satisfies $2dV=\sum_e Im(l_e)d\alpha_e$
(resp. $2dV=\sum_e l_ed\alpha_e$).
\end{lemma}

Note that the lengths considered here are those defined above, so that
they are in $(0,\pi)$ for spherical simplices, and in $\R_{<0}$ for
hyperbolic simplices.

\subsection{The generalized volume} \label{ssc:volume}

The previous lemma leads to a fairly natural definition of a real-valued
volume over the space of angled simplices.

\begin{defi}
Let $s$ be a tetrahedron and let $\omega$ be the $1$-form
(Schl\"afli 1-form) defined on $AS(s)$ by $2\omega=\sum_e
(Re(l_e)+Im(l_e)) d\alpha_e$.
\end{defi}

Note that the Schl\"afli 1-form is a continuous 1-form defined on the
6-dimensional convex polytope $AS(s)$.

It is proved in \cite{luo-ajm} that,

\begin{lemma}
$\omega$ is closed.
\end{lemma}

Remark that $\omega$ vanishes on the subspace of Euclidean
simplices.

\begin{defi}
The generalized volume $V:AS(s)\rightarrow \R$ is the primitive of
$\omega$ which vanishes on the Euclidean simplices.
\end{defi}

There is another possibility, namely to define the volume as a
complex-valued function, defining $\omega$ as $(1/2)\sum_e l_e d\theta_e$.
The definition chosen here serves well for our purposes.

Note that $V$ corresponds to the usual volume on
spherical and hyperbolic simplices by Lemma \ref{lm:schlafli}.
The volume of Euclidean simplices is zero by definition.
However, the volume of
flipped and doubly flipped Euclidean simplexes are not zero.

\begin{lemma} \label{lm:vol-euclidean}
Let $\alpha\in AS(s)$.
\begin{enumerate}
\item Suppose that $\alpha$ is flipped Euclidean, with the lengths
of the edges adjacent to $v_1$ equal to $i\pi$ and the other lengths
equal to $0$. Then
$$ V(\alpha)=\pi(\alpha_{12}+\alpha_{13}+\alpha_{14}-\pi)~. $$
\item Suppose that $\alpha$ is doubly flipped Euclidean, with
$l_{12}=l_{34}=0$ and the other lengths equal to $i\pi$. Then
$$ V(\alpha) = \pi(\alpha_{13}+\alpha_{14}+\alpha_{23}+\alpha_{24}-2\pi)~. $$
\end{enumerate}
\end{lemma}

Note that in each case the volume, without the factor $\pi$, is equal to
the area of a spherical polygon -- this will be useful below.

\begin{proof}[Proof of Lemma \ref{lm:vol-euclidean}]
For the first case, consider a small deformation that increases
slightly the $\alpha_{1i}$, $2\leq i\leq 4$. This deforms $\alpha$
into a spherical simplex $\alpha'$,  with vertices $v_2,v_3$ and $v_4$ very close
to the antipode of $v_1$. The (spherical) Schl\"afli formula, applied
to a 1-parameter formula deforming this simplex to a segment of length
$\pi$, shows that the volume of this simplex is equal to
$\pi(\alpha'_{12}+\alpha'_{13}+\alpha'_{14}-\pi)$, and the result follows
for $\alpha$.

The same argument works in the second case, the corresponding spherical
simplex now has $v_1,v_2$ very close and almost antipodal to both $v_3$
and $v_4$.
\end{proof}

There is a quite different way to define this ``volume'' of
domains in the extended hyperbolic space, in terms of an analytic
continuation \cite{cho-kim}.

\subsection{Smoothness}

For a closed triangulated 3-manifold $(M, T)$, the volume $V$ of
an angle system $x \in AS(T)$ is the sum of the volume of its
angled 3-simplexes. Thus $v: AS(T) \to \R$ is a $C^1$ smooth
function. Moreover it is real analytic outside the set of Euclidean
type simplices.

\section{Critical points}

This section contains the proof of Theorem \ref{tm:main}.

\subsection{Gluing conditions}

Suppose $(M, T)$ is a connected triangulated closed 3-manifold
so that $AS(T) \neq \emptyset$.  We will consider the volume
optimization $V: AS(T) \to \R$.

\begin{lemma}
Let $\theta_0\in AS(T)$ be a critical point of $V$ on $AS(T)$. Then, for each edge $e$ of
$T$, the lengths of $e$ for all the simplices containing it are equal.
\end{lemma}

This follows from the definition of $\omega$.  By the
classification lemma \ref{lm:simplex} and the connectivity of $M$, we see that
the types of any two 3-simplexes in $T$ in $\omega$ are the same.

If $\omega$ is a local maximum point of $V$, then it cannot happen
that all 3-simplices in $\omega$ are Euclidean simplices. Indeed,
if otherwise, the volume of $\omega$ is zero. However, if we
perturb $\omega$ slightly in $AS(T)$ to obtain a new point
$\omega'$, then all simplices in $\omega'$ can be hyperbolic and
spherical simplices. Thus $V(\omega') > V(\omega)$ which
contradicts the local maximum condition.

According to Lemma \ref{lm:lengths-unique}, for non-Euclidean type simplices, edge
lengths determine the isometry type. So we obtain,

\begin{cor}
Suppose $\omega$ is not of Euclidean type.
The faces of the simplices can be glued isometrically.
Furthermore,  $\theta_0$ defines in this way either a spherical
structure or a HS-structure $h$ on $M$.
\end{cor}

Indeed, there are two possibilities. Namely either all simplieces
in $\omega$ are of spherical type, or they are all of hyperbolic type. In the
spherical type case, all simplices are spherical and are glued by
isometries so that the sum of angles around each edge is $2\pi$.
Thus we obtain a spherical metric on $M$. In the case where all
simplices are of hyperbolic type, by Lemma \ref{lm:simplex}, we realize each
simplex in $\omega$ as a geometric tetrahedron in $\HSt^3$
so that their faces can be glued isometrically. Thus, we obtain an
HS-structure on $M$. Indeed, there are two subcases which could
occur. In the first case, all simplices are hyperbolic. Thus we
obtain a hyperbolic metric on $M$. In the second case, some
simplex is a flipped hyperbolic. Then we obtain an HS-structure on
$M$ by gluing these geometric tetrahedra in HS-geometry.

Note that all vertices of $T$ are in the hyperbolic part of this HS-structure.

\subsection{Normal spheres in HS-structures} \label{ssc:normal}

Continuing the proof of Theorem \ref{tm:main}, we consider here an HS-structure $h$ on $M$,
along with a triangulation $T$ with all vertices of $T$ in the hyperbolic part
of $h$. Suppose moreover that the de Sitter part $M_{dS}$ for $h$ is non-empty.
Let $M_0$ be a connected component of $M_{dS}$.

Then $M_0$ is geodesically complete, so it is isometric either to the de Sitter space
$dS^3$ or to its quotient by the antipodal map (see \cite{mess,mess-notes}).
Therefore any plane in the tangent space to $M_0$ at a point is tangent to
a (unique) totally geodesic space-like plane in $M_0$, which is homeomorphic
either to $S^2$ (in the first case) or to $\R P^2$ (in the second case).
Each of those totally geodesic surfaces is a normal surface in the triangulation $T$
of $M$.

This simple argument shows that each connected component of the de Sitter part
of $h$ corresponds to a normal surface in $(M,T)$.

\subsection{Normal spheres  for Euclidean critical points of $V$}

In this section we consider the same question as in \S \ref{ssc:normal}, about normal
surfaces in $(M,T)$, but for critical points of $V$ for which all simplices are of
Euclidean type. The arguments are of the somehow similar but are less geometric and more
combinatorial, because the geometric structures on the simplices cannot be glued to
obtain a geometric structure of Euclidean type on $M$.

We have seen in \S \ref{ssc:volume} that to each flipped (resp. doubly flipped)
Euclidean simplex $s$ in $T$ can be associated a spherical triangle (resp. a quadrilateral).
The edges of this triangle (resp. quadrilateral) are associated to the 2-faces
of $T$ which have exactly two edges of length $i\pi$. Each such face bounds
two simplices which are both either flipped or doubly flipped. It follows that
the triangles (resp. quadrilaterals) can be glued along their edges to obtain a
closed surface $\Sigma$ (which in general is not connected) -- however this gluing cannot
in general be isometric for the spherical metrics since the lengths of the edges
do not match. Moreover, the vertices of the triangulation of $\Sigma$ correspond
to the edges of $T$ of length $i\pi$.

\begin{remark}
  The angles of the triangles (resp.
quadrilateral) at each vertex sum up to $2\pi$.
\end{remark}

\begin{proof}
The angles of the triangles (resp. quadrilateral) adjacent to each vertex of $\Sigma$
are equal to the angles of the simplices of $T$ at the corresponding edge of length
$i\pi$. Those angle sum up to $2\pi$ by definition of a angle structure on $T$.
\end{proof}

\begin{cor}
Each connected component of $\Sigma$ is homeomorphic either to the sphere or to the
projective plane. The sum of the areas of the faces of $\Sigma$ is an integer
multiple of $2\pi$.
\end{cor}

\begin{proof}
Let $\Sigma_0$ be a connected component of $\Sigma$, let $F_0$ be the set of its
2-faces, and let $V_0$ be the set of its vertices. Given $f\in F_0$ and $v\in V_0$,
we write $v\simeq f$ if $v$ is adjacent to $f$, in this case we call $\theta_{f,v}$
the angle of $f$ at $v$.

Let $a(f)$ be the area of the face $f$ of $\Sigma$.
For each face $f\in F_0$ of $\Sigma_0$, we have by the Gauss-Bonnet formula
$$ \sum_{v\in V_0, v\simeq f} (\pi-\theta_{f,v}) = 2\pi - a(f)~. $$
Summing over the faces of $\Sigma_0$ yields that
$$ \sum_{f\in F_0}\left( \sum_{v\in V_0, v\simeq f} (\pi-\theta_{f,v})\right)
= 2\pi \# F_0 - \sum_{f\in F_0} a(f)~. $$
The number of wedges in the triangulation of $\Sigma_0$ is twice the number of
edges, which we denote by $\# E_0$. Therefore
$$ \sum_{f\in F_0}\sum_{v\in V_0, v\simeq f} \pi = 2\pi \# E_0~. $$
Moreover the angles of the faces at each vertex sum up to $2\pi$, so that
$$ \sum_{f\in F_0}\sum_{v\in V_0, v\simeq f} \theta_{f,v}) = 2\pi \# V_0~. $$
Using the definition of the Euler characteristic, we obtain that
$$ \sum_{f\in F_0} a(v) = 2\pi \# V_0 - 2\pi E_0 = 2\pi F_0 = 2\pi \chi(\Sigma_0)~, $$
and both parts of the corollary follow immediately.
\end{proof}

\begin{cor}
At a critical point of $V$ where all simplices are of Euclidean type, $V$ is
an integer multiple of $2\pi^2$.
\end{cor}

\begin{proof}
Lemma \ref{lm:vol-euclidean} shows that the volume of each flipped (resp.
doubly flipped) simplex is equal to $\pi$ times the area of the corresponding
triangle (resp. quadrilateral) in $\Sigma$. So the total volume is $\pi$ times
the area of $\Sigma$, so that it is a non-negative integer multiple of $2\pi^2$.
\end{proof}

The proof of Theorem \ref{tm:main} is obtained by putting together the results
of this section.

\section{Further questions}

The main point presented here is that extended hyperbolic
structures have a natural role when constructing geometric
structures on manifolds by maximization of volume over
triangulated manifolds. This leads to a number of questions, for
which answers would presumably help make progress on the
understanding of geometric structures on 3-manifolds.

\begin{q}
If M is a connected sum of several hyperbolic
3-manifolds, does M support an HS-structure?
\end{q}

Another, more general question, is whether the constructions considered here
can be extended to encompass angles structures with some ideal vertices. This
would mean allowing angle structures on simplices for which the sum of the
angles at a vertex is equal to, rather than less than, $2\pi$. Our hope is
that such ideal vertices would permit critical points of the volume to
realize torus decompositions of non atoroidal 3-manifolds. Another possibility,
adding some flexibility to the construction, would be to allow for vertices in
the de Sitter part of the extended hyperbolic space.

Another natural question is of course to understand the critical points of
$V$ on the boundary of $AS(T)$, hopefully showing that those boundary
critical points correspond to collapsings.

A last, more technical question, is whether existence of a
critical point of $V$ on $AS(T)$ for which all simplices are of
Euclidean type has topological consequences on $M$. For instance,
if all simplices are Euclidean (rather than only of Euclidean
type), does it follow that $M$ admits an Euclidean metric or more
generally is $M$ a connected sum of Seifert fibered spaces? This
is not obvious since the angles of the simplices add up to $2\pi$
at the edges of $T$, but the edge lengths do not match so that the
faces of the simplices cannot be isometrically glued.


\newcommand{\etalchar}[1]{$^{#1}$}
\def\cprime{$'$}

\end{document}